\documentclass[leqno,draft]{article}

\def\diam{\mathop{\rm diam}}

\begin{document}

\title{An introduction to the geometry \\ of metric spaces}

\author{Stephen Semmes \\
        Rice University}

\date{}

\maketitle

\begin{abstract}
These informal notes deal with some basic properties of metric spaces,
especially concerning lengths of curves.
\end{abstract}

\tableofcontents

\section{Definitions and notation}
\label{definitions, notation}
\setcounter{equation}{0}

        A \emph{metric space} is a nonempty set $M$ with a distance
function $d(x, y)$ defined for every $x, y \in M$.  More precisely,
$d(x, y)$ is supposed to be a nonnegative real number which is equal
to $0$ if and only if $x = y$, which is symmetric in $x$ and $y$ in
the sense that
\begin{equation}
        d(y, x) = d(x, y),
\end{equation}
and which satisfies the \emph{triangle inequality}
\begin{equation}
        d(x, z) \le d(x, y) + d(y, z)
\end{equation}
for every $x, y, z \in M$.  For example, the \emph{discrete metric} is
defined by putting $d(x, y) = 1$ when $x \ne y$, and it is easy to see
that this satisfies the preceding conditions.

        A more interesting example is the real line ${\bf R}$ with the
standard metric.  If $x$ is a real number, remember that the
\emph{absolute value} $|x|$ is defined by $|x| = x$ when $x \ge 0$ and
$|x| = - x$ when $x \le 0$.  It is well known and easy to check that
\begin{equation}
        |x + y| \le |x| + |y|
\end{equation}
for every $x, y \in {\bf R}$, and the standard metric on ${\bf R}$ is
defined by
\begin{equation}
        d(x, y) = |x - y|.
\end{equation}
If $(M, d(x, y))$ is any metric space and $E$ is a nonempty subset of
$M$, then the restriction of $d(x, y)$ to $x, y \in E$ defines a
metric on $E$, so that $E$ becomes a metric space too.

        Let $(M, d(x, y))$ be a metric space.  The \emph{open ball}
$B(x, r)$ with center $x \in M$ and radius $r > 0$ is defined by
\begin{equation}
        B(x, r) = \{y \in M : d(x, y) < r\}.
\end{equation}
Similarly, the \emph{closed ball} $\overline{B}(x, r)$ with center $x
\in M$ and radius $r \ge 0$ is defined by
\begin{equation}
        \overline{B}(x, r) = \{y \in M : d(x, y) \le r\}.
\end{equation}
Thus $\overline{B}(x, r)$ contains only $x$ when $r = 0$.

        A set $E \subseteq M$ is said to be \emph{bounded} if there is
a point $p \in M$ and a nonnegative real number $t$ such that
\begin{equation}
        d(p, x) \le t
\end{equation}
for every $x \in E$.  By the triangle inequality, this implies that
\begin{equation}
        d(q, x) \le d(p, q) + t
\end{equation}
for any $q \in M$ and $x \in E$.  One can use this to check that the
union of finitely many bounded sets is bounded.

\section{Norms on ${\bf R}^n$}
\label{norms}
\setcounter{equation}{0}

        Fix a positive integer $n$, and let ${\bf R}^n$ be the space
of $n$-tuples $x = (x_1, \ldots, x_n)$ of real numbers.  As usual, the
sum $x + y$ of $x, y \in {\bf R}^n$ is defined coordinatewise, so that
the $i$th coordinate of $x + y$ is the sum of the $i$th coordinates of
$x$ and $y$, $1 \le i \le n$.  If $t \in {\bf R}$ and $x \in {\bf
R}^n$, then the scalar product $t \, x$ is defined by putting its
$i$th coordinate equal to $t \, x_i$.

        A \emph{norm} on ${\bf R}^n$ is a nonnegative real-valued
function $N(x)$ defined for $x \in {\bf R}^n$ such that $N(x) = 0$ if
and only if $x = 0$,
\begin{equation}
        N(x + y) \le N(x) + N(y)
\end{equation}
for every $x, y \in {\bf R}^n$, and
\begin{equation}
        N(t \, x) = |t| \, N(x)
\end{equation}
for every $t \in {\bf R}$ and $x \in {\bf R}^n$.  In this case,
\begin{equation}
        d(x, y) = N(x - y)
\end{equation}
defines a metric on ${\bf R}^n$.

        For example, the absolute value function $|x|$ is a norm on
the real line.  The standard Euclidean norm on ${\bf R}^n$ is defined
by
\begin{equation}
        |x| = \Big(\sum_{j = 1}^n x_j^2\Big)^{1/2}.
\end{equation}
It is well known that this satisfies the triangle inequality, and
hence is a norm.  The corresponding metric is the standard Euclidean
metric on ${\bf R}^n$.  As another example,
\begin{equation}
        \|x\|_1 = \sum_{j = 1}^n |x_j|
\end{equation}
is a norm on ${\bf R}^n$.  For any real number $p \ge 1$, it can be shown that
\begin{equation}
        \|x\|_p = \Big(\sum_{j = 1}^n |x_j|^p \Big)^{1/p}
\end{equation}
is a norm on ${\bf R}^n$.  This is the same as the standard Euclidean
norm $|x|$ when $p = 2$, and the triangle inequality can be
established using the convexity of the function $|r|^p$ on ${\bf R}$
when $p \ge 1$.  It is easy to check directly that
\begin{equation}
        \|x\|_\infty = \max(|x_1|, \ldots, |x_n|)
\end{equation}
is a norm on ${\bf R}^n$.  Because
\begin{equation}
        \|x\|_\infty \le \|x\|_p
\end{equation}
and
\begin{equation}
        \|x\|_p \le n^{1/p} \, \|x\|_\infty,
\end{equation}
when $1 \le p < \infty$, $\|x\|_p \to \|x\|_\infty$ as $p \to \infty$.

         If $N$ is any norm on ${\bf R}^n$, then
\begin{equation}
         N(x) \le N(y) + N(x - y)
\end{equation}
and
\begin{equation}
         N(y) \le N(x) + N(x - y),
\end{equation}
and hence
\begin{equation}
         |N(x) - N(y)| \le N(x - y)
\end{equation}
for every $x, y \in {\bf R}^n$.  This implies that $N$ is continuous
with respect to the metric associated to $N$.  One can also check that
$N$ is bounded by a constant times the Euclidean norm on ${\bf R}^n$.
It follows that $N$ is continuous with respect to the standard
Euclidean metric on ${\bf R}^n$.

\section{The unit circle}
\label{unit circle}
\setcounter{equation}{0}

        The unit circle ${\bf S}^1$ is the set of $x = (x_1, x_2) \in
{\bf R}^2$ such that $|x| = 1$, or
\begin{equation}
        x_1^2 + x_2^2 = 1.
\end{equation}
The restriction of the standard Euclidean metric $|x - y|$ on ${\bf
R}^2$ to $x, y \in {\bf S}^1$ defines a metric on ${\bf S}^1$, but
there is another metric that is more intrinsic.  Specifically, let
$d(x, y)$ be the length of the shorter arc connecting $x$ to $y$ in
${\bf S}^1$.  This is the same as the angle at the origin between the
line segments to $x$ and $y$.  Clearly $d(x, y) \ge 0$ is symmetric in
$x$ and $y$, and is equal to $0$ exactly when $x = y$.  The total
length of the unit circle is $2 \pi$, and hence
\begin{equation}
        d(x, y) \le \pi
\end{equation}
for every $x, y \in {\bf S}^1$.  Furthermore,
\begin{equation}
        d(x, y) = \pi
\end{equation}
if and only if $x$ and $y$ are antipodal points in the circle, which
means that $y = -x$.  If $x, y, z \in {\bf S}^1$, then the shorter
arcs connecting $x$ to $y$ and $y$ to $z$ can be combined to get an
arc between $x$ and $z$, which implies that the triangle inequality
holds.

        There is a simple relationship between $d(x, y)$ and $|x - y|$,
which is that
\begin{equation}
        \sin \Big(\frac{d(x, y)}{2}\Big) = \frac{|x - y|}{2}.
\end{equation}
In particular,
\begin{equation}
        |x - y| \le d(x, y) \le \frac{\pi}{2} \, |x - y|
\end{equation}
for every $x, y \in {\bf S}^1$.  This can be improved when $|x - y|$
is small, since
\begin{equation}
        \lim_{r \to 0} \frac{\sin r}{r} = 1.
\end{equation}
Specifically, for every $\epsilon > 0$ there is a $\delta > 0$ such
that
\begin{equation}
        d(x, y) \le (1 + \epsilon) \, |x - y|
\end{equation}
when $|x - y| < \delta$.

        Put
\begin{equation}
        e(t) = (\cos t, \sin t)
\end{equation}
for each $t \in {\bf R}$.  This defines a mapping from the real line
onto the unit circle that satisfies
\begin{equation}
        e'(t) = \frac{d}{dt} e(t) = (\sin t, - \cos t).
\end{equation}
Hence
\begin{equation}
        |e'(t)| = 1
\end{equation}
for every $t \in {\bf R}$, which means that the length of the arc
traced by $e(t)$ for $a \le t \le b$ is equal to $b - a$ for every $a,
b \in {\bf R}$ with $a \le b$.  Thus
\begin{equation}
        d(e(a), e(b)) = b - a
\end{equation}
when $b - a \le \pi$.

\section{A little geometry}
\label{some geometry}
\setcounter{equation}{0}

        Let $P$ be a plane in ${\bf R}^n$.  If $x \in {\bf R}^n$ and
$x \not\in P$, then there is a unique point $x' \in P$ such that the
line in ${\bf R}^n$ passing through $x$ and $x'$ is perpendicular to
$P$.  If $x \in P$, then put $x' = x$.

        For each $p \in P$,
\begin{equation}
        |x - p|^2 = |x - x'|^2 + |x' - p|^2.
\end{equation}
If $r \ge |x - x'|$, then
\begin{equation}
        |x - p| = r
\end{equation}
is equivalent to
\begin{equation}
        |x' - p| = \widetilde{r}, \quad \widetilde{r}^2 = r^2 - |x - x'|^2
\end{equation}
for $p \in P$.

        Fix $q \in P$ and $t > 0$, and consider
\begin{equation}
        \{p \in P : |p - q| = t\}.
\end{equation}
The maximum and minimum of $|x' - p|$ on this set occur on the line
through $x'$ and $q$ when $x' \ne q$.  The maximum and minimum of
$|x - p|$ on this set occur at the same points.

\section{The unit sphere}
\label{unit sphere}
\setcounter{equation}{0}

        For each positive integer $n$, let ${\bf S}^{n - 1}$ be the
unit sphere in ${\bf R}^n$, consisting of the $x \in {\bf R}^n$ such
that $|x| = 1$.  This is the unit circle when $n = 2$, and it contains
only the two elements $\pm 1$ when $n = 1$.  Let us focus now on the
case where $n \ge 3$.  As before, the restriction of the Euclidean
metric $|x - y|$ to $x, y \in {\bf S}^{n - 1}$ defines a metric on
${\bf S}^{n - 1}$.  The \emph{spherical metric} $d(x, y)$ can be
defined be the conditions
\begin{equation}
        0 \le d(x, y) \le \pi
\end{equation}
and
\begin{equation}
        \sin \Big(\frac{d(x,y)}{2}\Big) = \frac{|x - y|}{2}
\end{equation}
for every $x, y \in {\bf S}^{n - 1}$.  This is symmetric in $x$ and
$y$, and equal to $0$ exactly when $x = y$.  In order to show that the
triangle inequality holds, we would like to reduce to the case of the
unit circle.

        Let $Q$ be a two-dimensional plane in ${\bf R}^n$ passing
through the origin.  The intersection of $Q$ with ${\bf S}^{n - 1}$ is
a circle of radius $1$, also known as a \emph{great circle} in ${\bf
S}^{n - 1}$.  We can think of $Q \cap {\bf S}^{n - 1}$ as a copy of
the unit circle, and $d(x, y)$ for $x, y \in Q \cap {\bf S}^{n - 1}$
corresponds exactly to the metric defined previously on ${\bf S}^1$.
Thus $d(x, y)$ also satisfies the triangle inequality on $Q \cap {\bf
S}^{n - 1}$.

        Let $x, y, z \in {\bf S}^{n - 1}$ be given, and let us show
that
\begin{equation}
        d(x, z) \le d(x, y) + d(y, z).
\end{equation}
This is trivial when $x = y$, and so we may suppose that $x \ne y$.
Let $Q$ be the two-dimensional plane in ${\bf R}^n$ passing through
$x$, $y$, and the origin.  Let $A$ be the set of $w \in {\bf S}^{n -
1}$ such that $|y - w| = |y - z|$, which is equivalent to $d(y, w) =
d(y, z)$.  This is an $(n - 2)$-dimensional sphere contained in ${\bf
S}^{n - 1}$, except for the trivial cases where $z = y$ or $z = -y$
and $A$ contains only $z$.  We can also describe $A$ as the
intersection of ${\bf S}^{n - 1}$ with a certain hyperplane $H$ in
${\bf R}^n$ perpendicular to the line $L$ through $y$ and $0$.  Using
geometric arguments as in the previous section, one can show that $|x
- w|$ is maximized on $A$ at a point $w_0 \in A \cap Q$.  Hence $d(x,
w)$ is maximized at the same point $w_0$.  We also have that
\begin{equation}
        d(x, w_0) \le d(x, y) + d(y, w_0),
\end{equation}
because $x$, $y$, and $w_0$ are contained in the same great circle $Q
\cap {\bf S}^{n - 1}$.  This implies that the triangle inequality
holds for $x$, $y$, and $z$, since $d(y, z) = d(y, w_0)$ by definition
of $A$ and $d(x, z) \le d(x, w_0)$ by maximization.

        As in the case of the unit circle,
\begin{equation}
        d(x, y) = \pi
\end{equation}
exactly when $x$ and $y$ are antipodal points in ${\bf S}^{n - 1}$,
which means that $y = -x$.  This is also equivalent to saying that $x$
and $y$ are contained in the same line passing through $0$.  For every
$x, y \in {\bf S}^{n - 1}$,
\begin{equation}
        |x - y| \le d(x, y) \le \frac{\pi}{2} \, |x - y|,
\end{equation}
and $d(x, y)$ is approximately the same as $|x - y|$ when $x$ and $y$
are close together, in the sense that for every $\epsilon > 0$
there is a $\delta > 0$ such that
\begin{equation}
        d(x, y) \le (1 + \epsilon) \, |x - y|
\end{equation}
when $|x - y| < \delta$.

\section{Supremum and infimum}
\label{supremum, infimum}
\setcounter{equation}{0}

        A real number $b$ is said to be an \emph{upper bound} for a set
$A \subseteq {\bf R}$ if $a \le b$ for every $a \in A$, and a real number
$c$ is said to be a \emph{lower bound} for $A$ if $c \le a$ for every
$a \in A$.  Note that $A$ has both an upper and lower bound in ${\bf R}$ if
and only if $A$ is bounded with respect to the standard metric on ${\bf R}$.

        A real number $\alpha$ is said to be the \emph{least upper
bound} or \emph{supremum} of a set $A \subseteq {\bf R}$ if $\alpha$
is an upper bound for $A$, and if $\alpha \le b$ for every upper bound
$b$ of $A$.  If $\alpha, \alpha' \in {\bf R}$ both satisfy these
conditions, then it follows that $\alpha \le \alpha'$ and $\alpha' \le
\alpha$, and hence $\alpha = \alpha'$.  The \emph{completeness
property} of the real numbers states that a nonempty set $A \subseteq
{\bf R}$ with an upper bound has a least upper bound, which is unique
by the previous remark.  The supremum of $A$ is denoted $\sup A$ when
it exists.  If $A$ has only finitely many elements, then the supremum
of $A$ is the same as the maximum of the elements of $A$.  Otherwise,
the supremum may not be an element of $A$.  For example, if $A$ is the
set of all negative real numbers, then $\sup A = 0$ is not an element
of $A$.

        The \emph{greatest lower bound} or \emph{infimum} of a set $A
\subseteq {\bf R}$ is defined analogously as a lower bound for $A$
which is greater than or equal to any other lower bound of $A$.  The
infimum of $A$ is unique when it exists, in which case it is denoted
$\inf A$.  It follows from the completeness property of the real
numbers that $A$ has an infimum when $A \ne \emptyset$ has a lower
bound.  Specifically, the infimum of $A$ can be obtained as the
supremum of the set of lower bounds for $A$.  Alternatively, the
infimum of $A$ is the negative of the supremum of $-A = \{-a : a \in
A\}$.

        A set $E$ in a metric space $(M, d(x, y))$ is bounded if and
only if the set of real numbers $d(x, y)$, $x, y \in E$, has an upper
bound.  If $E \subseteq M$ is bounded and nonempty, then the
\emph{diameter} $\diam E$ of $E$ is defined by
\begin{equation}
        \diam E = \sup \{d(x, y) : x, y \in E\}.
\end{equation}

\section{Lipschitz mappings}
\label{Lipschitz}
\setcounter{equation}{0}

        Let $(M, d(x, y))$ and $(N, \rho(u, v))$ be metric spaces.  A
mapping $f : M \to N$ is said to be \emph{Lipschitz} with constant $C
\ge 0$ if
\begin{equation}
        \rho(f(x), f(y)) \le C \, d(x, y)
\end{equation}
for every $x, y \in M$.  We may also simply say that $f$ is
$C$-Lipschitz in this case.  Thus a mapping is $0$-Lipschitz if and
only if it is constant.

        A mapping $f : M \to {\bf R}$ is $C$-Lipschitz with
respect to the standard metric on the real line if and only if
\begin{equation}
        f(x) \le f(y) + C \, d(x, y)
\end{equation}
for every $x, y \in M$.  This follows by interchanging the roles of $x$ and
$y$.  In particular, $f_p(x) = d(p, x)$ is $1$-Lipschitz for every $p \in M$.

        Lipschitz mappings are automatically uniformly continuous.  If
$f$ is a $C$-Lipschitz mapping from $M$ into $N$ and $E \subseteq M$
is nonempty and bounded, then
\begin{equation}
        {\diam}_N f(E) \le C \, {\diam}_M E,
\end{equation}
where the subscripts indicate in which metric space the diameter is taken.

        Let $(M_1, d_1)$, $(M_2, d_2)$, and $(M_3, d_3)$ be metric
spaces.  If $f_1 : M_1 \to M_2$ and $f_2 : M_2 \to M_3$ are Lipschitz
mappings with constants $C_1, C_2 \ge 0$, respectively, then the
composition $f_2 \circ f_1 : M_1 \to M_3$ defined by $(f_2 \circ
f_1)(x) = f_2(f_1(x))$ is Lipschitz with constant $C_1 \, C_2$.

\section{Lengths of curves}
\label{lengths, curves}
\setcounter{equation}{0}

        Let $(M, d(x, y))$ be a metric space, and let $a$, $b$ be real
numbers with $a \le b$.  The closed interval $[a, b]$ is defined as
usual as the set of real numbers $t$ such that $a \le t \le b$.  Let
$p : [a, b] \to M$ be a continuous mapping, which is to say a
continuous path in $M$ defined on $[a, b]$.

        A \emph{partition} $\mathcal{P}$ of $[a, b]$ is a finite
sequence $\{t_j\}_{j = 0}^n$ of real numbers such that
\begin{equation}
        a = t_0 < t_1 < \cdots < t_n = b.
\end{equation}
For each such partition $\mathcal{P}$, put
\begin{equation}
        \Lambda_\mathcal{P} = \sum_{j = 1}^n d(p(t_j), p(t_{j - 1})).
\end{equation}
This is an approximation to the length of $p(t)$, $a \le t \le b$.  If
$n = 1$, then
\begin{equation}
        \Lambda_\mathcal{P} = d(p(a), p(b)).
\end{equation}

        We say that $p$ has \emph{finite length} if the numbers
$\Lambda_\mathcal{P}$ have an upper bound, uniformly over all
partitions $\mathcal{P}$ of $[a, b]$.  In this case, the \emph{length}
of $p$ is denoted $\Lambda$ and defined to be the supremum of the
$\lambda_\mathcal{P}$'s.  Thus
\begin{equation}
        d(p(a), p(b)) \le \Lambda.
\end{equation}
Similarly, one can show that $p([a, b])$ is a bounded set in $M$ when
$p$ has finite length, and that
\begin{equation}
        \diam p([a, b]) \le \Lambda.
\end{equation}

         The condition of finite length is already nontrivial when $M$
is the real line equipped with the standard metric, for which it is
known classically as bounded variation.  The length of a real-valued
function is also known as the total variation.  If $p : [a, b] \to
{\bf R}$ is monotone increasing, then $p$ has bounded variation, and
the total variation of $p$ is equal to $p(b) - p(a)$.  However, one
can give examples of continuous real-valued functions on closed
intervals that do not have bounded variation.

\section{Special cases}
\label{cases}
\setcounter{equation}{0}

         Let $(M, d(x, y))$ be a metric space, and let $a$, $b$ be
real numbers with $a \le b$.  We can think of $[a, b]$ as being
equipped with the restriction of the standard metric on the real line.
If $p : [a, b] \to M$ is $C$-Lipschitz for some $C \ge 0$, then
$\Lambda_\mathcal{P} \le C \, (b - a)$ for every partition
$\mathcal{P}$ of $[a, b]$.  Thus $p$ has finite length $\Lambda \le C
\, (b - a)$.

         Suppose that $M$ is ${\bf R}^n$ with the standard Euclidean
metric, and that $p : [a, b] \to {\bf R}^n$ is continuously
differentiable.  The fundamental theorem of calculus implies that
\begin{equation}
         p(t) - p(r) = \int_r^t p'(u) \, du
\end{equation}
when $a \le r \le t \le b$.  Because $p'(u)$ is continuous on $[a,
b]$, $|p'(u)|$ is bounded on $[a, b]$, and $p$ is Lipschitz on $[a,
b]$.  It follows that $p$ has finite length on $[a, b]$.

         In this case, the length $\Lambda$ of $p$ on $[a, b]$ is
given by
\begin{equation}
         \Lambda = \int_a^b |p'(u)| \, du.
\end{equation}
For if $\mathcal{P}$ is any partition of $[a, b]$, then the previous
formula implies that
\begin{equation}
         \Lambda_\mathcal{P} \le \int_a^b |p'(u)| \, du.
\end{equation}
Hence
\begin{equation}
         \Lambda \le \int_a^b |p'(u)| \, du.
\end{equation}
To get the oppposite inequality, one can use uniform continuity of
$p'$ on $[a, b]$ to approximate the Riemann sums of the integral by
$\Lambda_\mathcal{P}$'s.

          There are analogous statements for the metric $d_N$
associated to a norm $N$ on ${\bf R}^n$.  Any norm on ${\bf R}^n$ is
bounded by a constant multiple of the standard norm, which implies
that a continuously-differentiable curve $p : [a, b] \to {\bf R}^n$ is
also Lipschitz with respect to the metric $d_N$ on ${\bf R}^n$, and
thus has finite length.  The norm of the integral of a continuous
${\bf R}^n$-valued function is less than or equal to the integral of
the norm of the function, as in the case of the Euclidean norm, and hence
\begin{equation}
          N(p(r) - p(t)) \le \int_r^t N(p'(u)) \, du
\end{equation}
when $a \le r \le t \le b$.  This implies that
\begin{equation}
          \Lambda_\mathcal{P} \le \int_a^b N(p'(u)) \, du
\end{equation}
for any partition $\mathcal{P}$ of $[a, b]$, where
$\Lambda_\mathcal{P}$ is now the approximation to the length of the
curve corresponding to the metric $d_N$.  Therefore
\begin{equation}
          \Lambda \le \int_a^b N(p'(u)) \, du,
\end{equation}
and one can get the opposite inequality to conclude that
\begin{equation}
          \Lambda = \int_a^b N(p'(u)) \, du
\end{equation}
in the same way as for the Euclidean norm.

\section{Compositions}
\label{compositions}
\setcounter{equation}{0}

        Let $(M, d(x, y))$ and $(N, \rho(u, v))$ be metric spaces,
and suppose that $f$ is a Lipschitz mapping from $M$ to $N$ with
constant $C \ge 0$.  If $p$ is a continuous mapping from $[a, b]$ into
$M$, then the composition $f \circ p$ is a continuous mapping from
$[a, b]$ into $M$.  For each partition $\mathcal{P}$ of $[a, b]$, the
analogue of $\Lambda_\mathcal{P}$ for $f \circ p$ is bounded by $C$
times $\Lambda_\mathcal{P}$ for $p$.  It follows that $f \circ p$ has
finite length if $p$ has finite length, and that the length of $f
\circ p$ is bounded by $C$ times the length of $p$.

        Suppose that $\alpha$, $\beta$ are also real numbers such
that $\alpha \le \beta$, and that $\phi$ is a one-to-one continuous
mapping of $[\alpha, \beta]$ onto $[a, b]$.  Thus $\phi$ maps every
partition $\mathcal{P}$ of $[\alpha, \beta]$ to a partition
$\phi(\mathcal{P})$ of $[a, b]$, and every partition of $[a, b]$ is of
the form $\phi(\mathcal{P})$ for some partition $\mathcal{P}$ of
$[\alpha, \beta]$.  The approximation to the length of $p \circ \phi$
associated to a partition $\mathcal{P}$ of $[\alpha, \beta]$ is equal
to the approximation to the length of $p$ associated to the
corresponding partition $\phi(\mathcal{P})$ of $[a, b]$.  It follows
that $p : [a, b] \to M$ has finite length if and only if $p \circ \phi
: [\alpha, \beta] \to M$ does, in which event the lengths are the
same.

        Now suppose that $\phi$ is a monotone increasing continuous
mapping from $[\alpha, \beta]$ onto $[a, b]$.  If $\alpha \le \rho \le
\tau \le \beta$ and $\phi(\rho) = \phi(\tau)$, then $\phi$ is constant
on $[\rho, \tau]$.  Although the correspondence between partitions of
$[\alpha, \beta]$ and $[a, b]$ is a little more complicated, the
approximations to the lengths of $p$ and $p \circ \phi$ still match
up, because the intervals on which $\phi$ is constant only add terms
equal to $0$ to the approximations to the length of $p \circ \phi$.
Consequently, $p$ has finite length if and only if $p \circ \phi$
does, and the lengths are again the same.

        One can also consider the composition $p \circ \phi$ when
$\phi$ is a continuous mapping from $[\alpha, \beta]$ onto $[a, b]$
which may not be monotone.  If $p \circ \phi$ has finite length, then
one can check that $p$ has finite length, and that the length of $p$
is less than or equal to the length of $p \circ \phi$.  It is easy to
have strict inequality, because $\phi$ may retrace parts of $[a, b]$
more than once.  For example, if $M$ is the real line and $p$ is the
identity mapping, then $\phi = p \circ \phi$ may have unbounded
variation or total variation strictly greater than $b - a$.

\section{Refinements and subintervals}
\label{refinements, subintervals}
\setcounter{equation}{0}

        Let $(M, d(x, y))$ be a metric space, and let $p$ be a
continuous mapping from a closed interval $[a, b]$ in the real line
into $M$.  A partition $\mathcal{P}_2$ of $[a, b]$ is said to be a
\emph{refinement} of another partition $\mathcal{P}_1$ of $[a, b]$ if
each point in $\mathcal{P}_1$ is also in $\mathcal{P}_2$.  If
$\Lambda_{\mathcal{P}_1}$, $\Lambda_{\mathcal{P}_2}$ are the
corresponding approximations to the length of $p$, then one can use
the triangle inequality to show that
\begin{equation}
        \Lambda_{\mathcal{P}_1} \le \Lambda_{\mathcal{P}_2}.
\end{equation}
Note that for every pair of partitions $\mathcal{P}$, $\mathcal{P}'$
of $[a, b]$, there is a partition $\mathcal{P}''$ of $[a, b]$ which is
a refinement of both $\mathcal{P}$ and $\mathcal{P}'$.

         If $a_1$, $b_1$ are real numbers such that $a \le a_1 \le b_1
\le b$, then every partition of $[a_1, b_1]$ can be extended to a
partition of $[a, b]$.  The approximation to the length of $p$ on
$[a_1, b_1]$ corresponding to the first partition is less than or
equal to the approximation to the length of $p$ on $[a, b]$ that
corresponds to the second partition.  If $p$ has finite length on $[a,
b]$, then it follows that the restriction of $p$ to $[a_1, b_1]$ has
finite length less than or equal to the length of $p$ on $[a, b]$.
Let the length of $p$ on $[a_1, b_1]$ be denoted $\Lambda(a_1, b_1)$,
so that the previous statement is expressed by the inequality
\begin{equation}
       \Lambda(a_1, b_1) \le \Lambda(a, b).
\end{equation}

        For each $r \in [a, b]$,
\begin{equation}
        \Lambda(a, r) + \Lambda(r, b) = \Lambda(a, b).
\end{equation}
Any partitions of $[a, r]$ and $[r, b]$ can be combined to get a
partition of $[a, b]$, and the sum of the corresponding approximations
to $\Lambda(a, r)$ and $\Lambda(r, b)$ is an approximation to
$\Lambda(a, b)$, which implies that $\Lambda(a, r) + \Lambda(r, b) \le
\Lambda(a, b)$.  Every partition of $[a, b]$ has a refinement of this
form, which implies the opposite inequality.  The same argument shows
that $p$ has finite length on $[a, b]$ when its restrictions to $[a,
r]$ and $[r, b]$ have finite length.

        As a function of $r$ on $[a, b]$, $\Lambda(a, r)$ is monotone
increasing, and one can also show that $\Lambda(a, r)$ is continuous.
It suffices to check that $\Lambda(a, r)$ is continuous from the left,
and that $\Lambda(r, b)$ is continuous from the right.  To do this, it
is helpful to consider partitions of $[a, r]$ and $[r, b]$ for which
the corresponding approximations to the length are close to the
supremum, and to use the previous remarks about refinements of
partitions.  The continuity of $p$ is also important here, to limit
the effect of a term in the sums involving $r$.

\section{Curves of minimal length}
\label{minimal}
\setcounter{equation}{0}

        Let $(M, d(x, y))$ be a metric space, and suppose that $p :
[a, b] \to M$ is a continuous curve of finite length.  If $\Lambda(a,
r)$ is the length of the restriction of $p$ to $[a, r]$, then there is
a continuous mapping $q : [0, \Lambda(a, b)] \to M$ such that
$q(\Lambda(a, r)) = p(r)$.  The main point is that $p(r)$ is constant
on any interval on which $\Lambda(a, r)$ is constant, so that $q$ is
well-defined.  Moreover, $q$ is Lipschitz with constant $1$, because
the distance between the endpoints of a curve is less than or equal to
the length of the curve.

        If closed and bounded subsets of $M$ are compact, then there
is a continuous curve of minimal length connecting any pair of
elements of $M$ for which there is a continuous curve of finite
length.  Using the observation in the preceding paragraph and a linear
change of variables on the real line, it is enough to show that there
is a Lipschitz mapping from the unit interval $[0, 1]$ into $M$
connecting the two points whose Lipschitz constant is as small as
possible.  One starts with a sequence of Lipschitz mappings from $[0,
1]$ into $M$ connecting the two points whose Lipschitz constants tend
to the infimum.  The Arzela-Ascoli theorem implies that a subsequence
of this sequence converges uniformly to a Lipschitz mapping on $[0,
1]$ with minimal Lipschitz constant.

        For example, line segments in ${\bf R}^n$ yield paths of
minimal length for the standard metric, or for the metric associated
to any norm.  In any metric space, the length of a path is greater
than or equal to the distance between the endpoints of the path.  If
the length of the path is equal to the distance between the endpoints,
then the path automatically has minimal length.  This is exactly what
happens for line segments in ${\bf R}^n$ with respect to any norm,
although for some norms like $\|\cdot\|_1$ and $\|\cdot\|_\infty$
there are other paths of minimal length too.

        A curve in the unit sphere in ${\bf R}^n$ has finite length
with respect to the Euclidean metric if and only if it has finite
length with respect to the spherical metric, in which case the length
of the curve is the same for both metrics.  This is because the two
metrics are approximately the same locally in such a precise way, and
because one can use refinements of partitions to make sure that only
local distances are used to determine the length of a path.  An arc of
a great circle in ${\bf S}^{n - 1}$ with length $\le \pi$ has minimal
length with respect to the spherical metric, because the distance
between its endpoints is equal to its length.  Such an arc therefore
has minimal length with respect to the Euclidean metric as well.

\end{document}